\documentclass[12pt]{article}
\usepackage{epsfig}
\begin{document}

\newcommand{\nc}{\newcommand}
\nc{\eqn}[1]{Eq.(\ref{#1})}
\nc{\wph}{\vphantom{A^A}}
\nc{\vph}{\vphantom{{A^A\over A}}}
\nc{\bq}{\begin{equation}}
\nc{\eq}{\end{equation}}
\nc{\bqa}{\begin{eqnarray}}
\nc{\eqa}{\end{eqnarray}}
\nc{\nl}{\nonumber \\}    
\nc{\suml}{\sum\limits}
\nc{\prodl}{\prod\limits}
\nc{\avg}[1]{\left\langle #1\right\rangle}
\nc{\vari}[1]{\sigma\left(#1\right)^2}
\nc{\val}[1]{N^{{\underline{#1}}}}
\nc{\vam}[2]{#1^{{\underline{#2}}}}
\nc{\ddel}[1]{\delta\left(#1\right)}
\nc{\ddelt}[2]{\delta^{#1}\left(#2\right)}
\nc{\plaat}[4]{\raisebox{#4cm}{\epsfig{figure=figures/#1.eps,width=#2cm,height=#3cm}}}
\nc{\boxnobox}[2]{\begin{minipage}{#1}{#2}\end{minipage}}
\nc{\intl}{\int\limits}
\nc{\vx}{\vec{x}}
\nc{\vp}{\vec{p}}
\nc{\vq}{\vec{q}}
\nc{\vr}{\vec{r}}
\nc{\vy}{\vec{y}}
\nc{\vz}{\vec{z}}
\nc{\vecv}{\vec{v}}
\nc{\vn}{\vec{n}}
\nc{\sn}{_{\vec{n}}}
\nc{\vnul}{\vec{0}}
\nc{\al}{\alpha}
\nc{\be}{\beta}
\nc{\si}{\sigma}
\nc{\del}{\delta}
\nc{\mod}{\;\mbox{mod}\;}
\nc{\bino}[2]{\left(\begin{tabular}{c} $#1$ \\ $#2$ \end{tabular}\right)}
\nc{\ord}[1]{{\cal O}\left(#1\right)}
\nc{\vur}{V_{\mbox{UR}}}
\nc{\vnr}{V_{\mbox{NR}}}
\nc{\from}{\leftarrow}
\nc{\mc}{Monte Carlo}
\nc{\qmc}{Quasi-Monte Carlo}
\nc{\alg}{algorithm}
\nc{\algs}{algorithms}
\nc{\diag}[3]{\raisebox{#3cm}{\epsfig{figure=#1.eps,width=#2cm}}}

\begin{center}
{\Large {\bf First- and second-order error estimates\\ in Monte Carlo integration}}\\
\vspace*{2\baselineskip}
R. Bakx\footnote{{\tt renz.bakx@student.ru.nl}}, 
R.H.P. Kleiss\footnote{{\tt R.Kleiss@science.ru.nl}}, 
F. Versteegen\footnote{{\tt f.versteegen@student.ru.nl}}\\
Radboud University Nijmegen,\\
Institute for Mathematics, Astrophysics and Particle Physics, \\Heyendaalseweg 135, 
NL-6525 AJ Nijmegen, The Netherlands.\\
\vspace*{3\baselineskip}
{\bf Abstract}\\
\vspace*{\baselineskip}
\begin{minipage}{10cm}{{\small 
In Monte Carlo integration an accurate and reliable determination of the numerical intregration
error is essential. We point out the need for an independent estimate of the error on this error,
for which we present an unbiased estimator. In contrast to the usual (first-order) error estimator, this
second-order estimator can be shown to be not necessarily positive in an actual Monte Carlo
computation. We propose an alternative and indicate how this can be computed in linear time
without risk of large rounding errors. In addition, we comment on the relatively very slow convergence of the second-order error estimate.
}}\end{minipage}
\vspace*{2\baselineskip}\end{center}

\section{Monte Carlo integration and its errors}
It does not need to be stressed that in numerical integration, including Monte Carlo (MC) 
integration \cite{MCoriginal},
a determination or estimate of the integration error made is essential. The Central Limit Theorem (CLT)
practically ensures that if the number $N$ of MC points is sufficiently large the numerical value of
the MC integral - itself a stochastic variable - will have a Gaussian distribution around the true integral
value, with a standard deviation that can itself also be estimated: this is the {\em first-order error}. The
results of MC integrations are therefore usually reported as
\[
\mbox{"result"}\;\;\pm\;\;\mbox{"error"}
\]
with the understanding that the "error" value quoted is the Gaussian's standard deviation. In this way
one can, for instance, assign confidence levels when comparing the integration result with a 
measurement. However, since the Gaussian distribution is quite steep, a modest change in the
value of the error can change the confidence levels considerably. It is therefore preferable to also have
a {\em second-order error} that estimates how well the first-order error was computed. The better way to 
report the result of a MC integration is then
\[
\mbox{"result"}\;\;\pm\;\;\left(\mbox{"first-order error"}\;\;\pm\;\;\mbox{"second-order error"}\vph\right)\;\;.
\]
A first attempt to implement such a method was presented in \cite{BKP}. However, in that paper no explicit form of the second-order error estimator was presented, nor
were its numerical stability properties and its convergence behaviour discussed:
also it was (wrongly) stated that the second-order error was the square root of the
estimator, while it ought to be the fourth root. The present paper addresses and
corrects these issues.
In what follows we shall arrive at an estimator for the second-order error that, like the first-order one,
can be evaluated in linear time {\it i.e.\/} at essentially no extra CPU cost. We shall also discuss
several of its numerical aspects, and suggest an improvement.

\section{Error estimators}
We will start by defining some mathematical tools. We consider an integral over an integration region
$\Gamma$ of an integrand $f(x)$, with $x\in\Gamma$. We have at our disposal a set of MC
integration points $x_j\;,j=1,2,\ldots,N$, 
assumed to be iid (Independent, Identically Distributed)
with a probability distribution $P(x)$ in $\Gamma$. We define
\bq
J_p = \intl_\Gamma dx\,P(x)\;w(x)^p\;\;\;,\;\;\;w(x) = {f(x)\over P(x)}\;\;,
\eq
so that $J_1 = \int dx \;f(x)$, the integral we want to compute. The numbers $w_j\equiv w(x_j)$ 
are called the
{\em weights\/} of the points. We see that $J_p$ is nothing but the expectation value of
$w(x)^p$:
\bq
\avg{w^p} = J_p\;\;.
\eq

Furthermore, we define the following multiple sums:
\bq
S_{p_1,p_2,\ldots,p_k} = \suml_{j_{1,2,\ldots,k}=1}^N {w_{j_1}}^{p_1}\; {w_{j_2}}^{p_2} \cdots
\cdots{w_{j_k}}^{p_k}
\eq
with the condition that the indices $j_{1,2,\ldots,k}$ are all {\em different}. As an example, the sum
$S_{1,1}$ does not contain $N^2$ but $\val{2} = N^2-N$ terms.
The {\em falling powers\/} are defined
by
\bq
\val{p} = N!/(N-p)! = N(N-1)(N-2)\cdots(N-p+1)\;\;.
\eq
The simple sums $S_p$ can be evaluated in linear time (that is, using $N$ additions), but a
multiple sum $S_{p_1,\ldots,p_k}$ needs time of the order $N^k$. In calculating estimators we 
therefore want to use only simple sums. On the other hand, only the multiple sums have a simple
expectation value:
\bq
\avg{S_{p_1,p_2,\ldots,p_k}} = \val{k}\;J_{p_1}\,J_{p_2}\cdots J_{p_k}\;\;.
\eq
We can relate simple and multiple sums to one another by the following obvious rule:
\bqa
S_{p_1,p_2,\ldots,p_k}S_q &=& S_{p_1+q,p_2,\ldots,p_k} + S_{p_1,p_2+q,\ldots,p_k}
+ \cdots + S_{p_1,p_2,\ldots,p_k+q}\nl && + \;S_{p_1,p_2,\ldots,p_k,q}\;\;.
\eqa

We are now ready to construct the various estimators, starting with the well-known MC 
formul\ae\ for clarity. For the integral we have
\bq
E_1 = {1\over N}S_1\;\;,
\label{e1estimator}
\eq
since $\avg{E_1}=J_1$; moreover we see that this estimator is unbiased. For the variance of $E_1$
we have
\bqa
\avg{{E_1}^2}-\avg{E_1}^2 &=& {1\over N^2}\avg{S_2 + S_{1,1}} - {J_1}^2\nl
&=& {1\over N}\left(J_2 - {J_1}^2\right)
= {1\over N^2}\avg{S_2} - {1\over\val{2}N}\avg{S_{1,1}}
\label{varianceofe1}
\eqa
so that the appropriate estimator is
\bq
E_2 = {S_2\over N^2} - {S_{1,1}\over\val{2}N} = {1\over\val{2}N}\Sigma_2\;\;\;,\;\;\;
\Sigma_2 = N\,S_2 - {S_1}^2\;\;.
\label{e2estimator}
\eq
The latter form is more suited to computation since it can be evaluated in linear time.
From \eqn{varianceofe1} we see that the first-order error, defined as
${E_2}^{1/2}$ decreases as $N^{-1/2}$, as is of course
very well known. Moreover, the expected error is defined for all
functions that are quadratically integrable, as is equally well known.\\

The {\em second-order error\/}
 should have as {\em its\/} expectation value the variance of $E_2$, which by the same methods as above can be shown to be
\bqa
\avg{{E_2}^2}-\avg{E_2}^2 &=& 
{1\over N^3}\left(J_4-4J_3J_1 + 3{J_2}^2 - 4\left(J_2-{J_1}^2\right)^2\wph\right)\nl &&
+ {2\over\val{2}N^2}\left(J_2-{J_1}^2\right)^2\;\;.
\eqa
We see that the second-order error, defined as ${E_4}^{1/4}$
 decreases, for large $N$, as $N^{-3/4}$. Moreover we see that the second-order error
 is only meaningful for integrands that are at least {\em quartically\/} integrable.
The appropriate unbiased estimator with the correct expectation value is 
\bqa
E_4 &=& {1\over\val{4}N^3}\left(\val{2}\Sigma_4 - 4{\Sigma_2}^2\right) + {2\over\val{4}\val{2}N^2}{\Sigma_2}^2\;\;,\nl
\Sigma_4 &=& N\,S_4 - 4\,S_3\,S_1 + 3\,{S_2}^2\;\;.
\label{e4estimator}
\eqa
An important observation here concerns the asymptotic behaviour of the
relative errors. Whereas the relative first-order error, {\it i.e.\/}
the ratio ${E_2}^{1/2}/E_1$, goes as $N^{-1/2}$ according to the `standard' behaviour in MC, the
relative second-order error ${E_4}^{1/4}/{E_2}^{1/2}$ only decreases as fast as $N^{-1/4}$. It will
therefore take much longer for the error to be well-determined than for the integral 
itself\footnote{Note that the relative errors as defined here are the {\em dimensionless\/}
 ratios, the
only meaningful measures of performance of the computation.}.

A final point is in order. By the CLT we know that the
distribution of $E_1$ in an ensemble of MC computations is normally distributed, 
which tells us the {\em meaning\/} of $E_2$, as discussed above. Since
$E_2$ is not computed as a simple average, {\em its\/} distribution is not governed by
the {\em same\/} CLT. Nevertheless, as is shown in the Appendix a good case can 
be made for it being also approximately normally distributed, so that the relation
between $E_4$ and the confidence levels of $E_2$ can be treated in the
usual manner. Below, we shall illustrate this with several examples.

\section{Positivity and numerical stability}
In principle, equations (\ref{e1estimator}), (\ref{e2estimator}) and (\ref{e4estimator}) are what
is necessary to obtain the integral and its first- and second-order errors. However, a number of
considerations must modify this picture. In the first place, the issue of positivity. 
Writing $w(x)=J_1 + u(x)$ so that $\int dx\,P(x)\,u(x) = 0$, we have
\bqa
\lefteqn{J_2 - {J_1}^2  = \int dx\;P(x)\,u(x)^2\;\;,}\nl
\lefteqn{J_4-4J_3J_1 + 3{J_2}^2 = \int dx\;P(x)\,u(x)^4 + 3\left(\int dx\;P(x)\,u(x)^2\right)^2\;\;,}\nl
\lefteqn{J_4-4J_3J_1 + 3{J_2}^2 - 4\left(J_2 - {J_1}^2\right)^2 =} \nl &&
{1\over2}\int dx\,dy\;P(x)\,P(y)\,(u(x)^2-u(y)^2)^2\;\;,\nl
\eqa
so that the {\em expectation values\/} of $E_{2,4}$ are positive, as they should. In
addition, since with the notation $W_j = E_1 + u_j$ the $\Sigma_2$ can be written as
\bq
\Sigma_2 = {1\over2}\suml_{j,k}\left(u_j-u_k\right)^2\;\;,
\eq
also $E_2$ itself is strictly nonnegative in any actual MC calculation. For $E_4$ this does
not hold, however. A counterexample can be constructed as follows. Let us assume that
the MC weights $w_j$ take on only the values 0 and 1, and that $E_1=Nb$, $b \in [0,1]$.
We then have
\bq
\Sigma_2 = \Sigma_4 = N^2a\;\;\;,\;\;\;a = b-b^2 \in [0,1/4]\;\;.
\eq
The value of $E_4$ now comes out as
\bq
E_4 = {1\over\val{4}}\left({\val{2}\over N}a - {4N^3-6N^2\over\val{2}}a^2\right)\;\;,
\eq
which is actually {\em negative\/} for
\bq
a >  {(N-1)^2\over N(4N-6)} = {1\over4} - {N-2\over2N(4N-6)}\;\;.
\eq
Although by small margin (surprisingly, in this counterexample, for $b\approx1/2$), the
positivity of $E_4$ cannot be guaranteed, so that ${E_4}^{1/4}$ may be undefined.
As an improvement on this situation we propose to abandon the estimator $E_4$ in favour of
\bq
\hat{E}_4 = {1\over\val{4}N^3}\left(\wph N^2\Sigma_4 - 4{\Sigma_2}^2\right)\;\;.
\eq
This estimator has a slight (order $1/N$) bias, which ought to be acceptable since
we are dealing with only the second-order error here; its advantage is that, since
\bq
N^2\Sigma_4 - 4{\Sigma_2}^2 =  {N^2\over2}\suml_{j,k}\left({u_j}^2-{u_k}^2\right)^2\;\;,
\eq
it always evaluates to a nonnegative number.\\

The second issue is that of numerical stability. It is well known that already the evaluation of
$\Sigma_2$ involves large cancellations which may destroy the numerical stability of the calculation
and can actually lead to negative values for $E_2$: this is the
reason why the straightforward computation of $E_2$ usually cannot be
reliably performed with single-precision arithmetic\footnote{As anyone who has
ever taught courses on Monte Carlo integration can testify.}. This problem has been 
widely discussed, for instance in \cite{James,CG}.
The situation of $E_4$, which involves even 
larger cancellations, is certainly worse. To tackle these problems, we adopt the CGV algorithm
first described in \cite{CG}. The strategy of this algorithm can best be summarized
as follows. In the first place, one concentrates on objects that are supposed to go to a 
finite asymptotic value. $E_1$ is such an object, but $\Sigma_{2,4}$ are not. 
In the second place, the algorithm
focuses on the update of these numbers as $N$ is increased by 1. So let us define
\bqa
M(N) &=& S_1(N)/N\;\;,\nl
P(N) &=& S_2(N)/N - S_1(N)^2/N^2\;\;,\nl
Q(N) &=& S_3(N)/N - 3S_2(N)S_1(N)/N^2 + 2S_1(N)^3/N^3\;\;,\nl
R(N) &=& S_4(N)/N - 4S_3(N)S_1(N)/N^2 + 3S_2(N)^2/N^2 - 4P(N)^2 .
\eqa
Here we have explicitly indicated the $N$ dependence of the running sums $S_{1,2,3,4}$.
We also define 
\bq
m = M(N-1)\;\;\;,\;\;\;p = P(N-1)\;\;\;,\;\;\;q = Q(N-1)\;\;\;,\;\;\;u = w_N - m\;\;.
\eq
The authors of \cite{CG} have already established the update rules
\bqa
M(N) = m + {1\over N}u\;\;,\nl
P(N) = {N-1\over N}\left( p + {1\over N}u^2\right)\;\;.
\eqa
We see that in particular the computation of $P(N)$ is free of large cancellations.
Some algebra leads us to supplement these update rules by
\bqa
Q(N) &=& {N-1\over N}\left( q + {N-2\over N^2}u^3  - {3p\over N}u\right)\;\;,\nl
R(N) &=& {N-1\over N}\left( R(N-1) + {1\over N}\left(p - {N-2\over N}u^2\right)^2
- 4\left({q\over N}u - {p\over N^2}u^2\right)\right)\;\;.\nl
\eqa
Using these results, for any given $N$ we then have
\bq
E_2 = {N\over\val{2}}P(N)\;\;\;,\;\;\;
\hat{E}_4 = {N\over\val{4}}R(N)\;\;.
\eq

\section{A case study}
To illustrate all the above, we can perform a simple but enlightening study.
Let us consider the following class of integrands:
\bq
f_\al(x) = (1+\al)\,x^\al\;\;\;,\;\;\;x \in (0,1]\;\;,\;\; -1 <\al \le 0\;\;,
\eq
which we shall integrate by employing $N$ pseudorandom numbers,  iid
uniformly in $(0,1]$. 
These functions are all integrable (with $J_1=1$), but divergent as $x\to 0$.
For $\al\le-0.5$ they are not quadratically integrable, and for $\al\le-0.25$ they
are quadratically integrable but not quartically integrable. Consequently, for
$\al\le-0.25$ the expectation value $\avg{\hat{E}_4}$ is not defined, and for
$\al\le-0.5$ not even $\avg{E_2}$ is defined. Nevertheless, in any actual MC
calculation of this integral, $\Sigma_{2,4}$ and $E_2,\hat{E}_4$ will have
definite, well-defined numerical values. So how, then, are these to be interpreted?

Below, we give the results for ${E_2}^{1/2}$ and ${\hat{E}_4}^{1/4}$ in a MC run
where $N\le 10^4$, monitoring their behaviour while $N$ increases. This we do for
values of $\al$ running from $-0.1$ down to $-0.9$. 
 
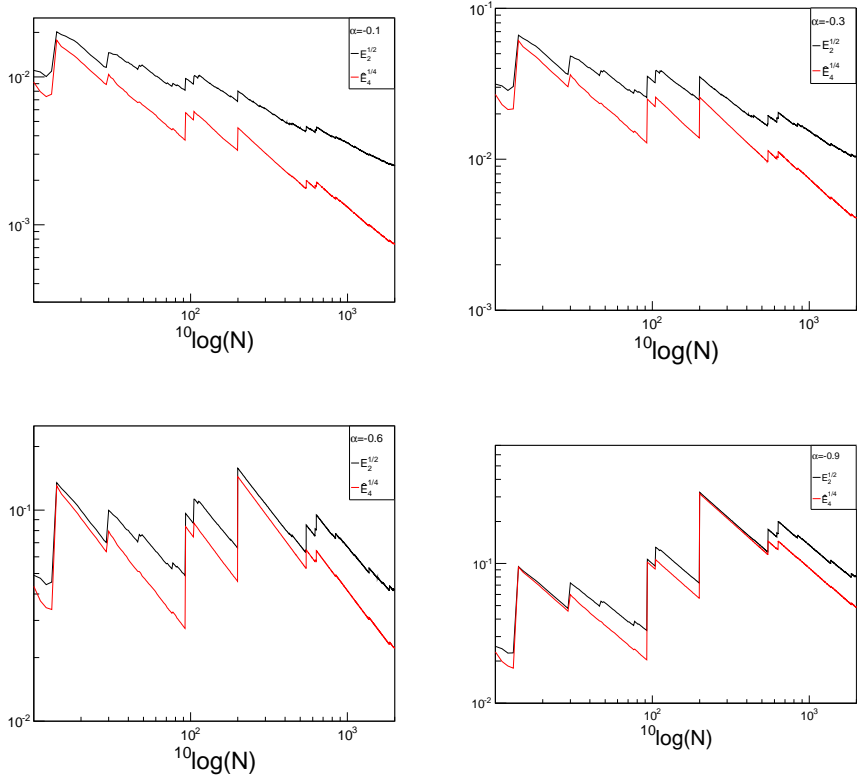
\begin{minipage}{6cm}{\[\diag{Plaatje_a01}{6}{0}\]}\end{minipage}
\begin{minipage}{6cm}{\[\diag{Plaatje_a03}{6}{0}\]}\end{minipage}

\begin{minipage}{6cm}{\[\diag{Plaatje_a06}{6}{0}\]}\end{minipage}
\begin{minipage}{6cm}{\[\diag{Plaatje_a09}{6}{0}\]}\end{minipage}

\noindent The upper line is the evolution of ${E_2}^{1/2}$, and the lower
line displays ${\hat{E}_4}^{1/4}$. For the smoothest case, $\al=-0.1$, the $N^{-1}$ behaviour for $E_2$ and the
$N^{-3}$ behaviour for $\hat{E}_4$ are evident\footnote{Note that in the plots
the values given are those of $E_2^{1/2}$ and $\hat{E}_4^{1/4}$.}, 
marred by smallish jumps whenever an $x$
value close to the singularity at $x=0$ is encountered. As $\al$ decreases to -0.3 quartic 
integrability is lost, which can be seen from the fact that the jumps in $\hat{E}_4$ are
now much larger while those in $E_2$ remain modest. Note that in all cases
exactly the same set of pseudrandom numbers was used. Therefore in the various
plots the jumps are in the same place, they simply become larger and larger. 
For $\al=-0.6$ where the integrand is also no 
longer quadratically integrable even the $N^{-1}$ behaviour of $E_2$ becomes quite distorted
by the growing jumps. Finally, at $\al=-0.9$ where the function itself is barely integrable,
the jumps have become so large that the short-term $N^{-1}$ and $N^{-3}$ behaviour inbetween
the jumps can
no longer ensure this behaviour over longer $N$ ranges. It is this kind of behaviour --- short-range
smooth decrease interspersed with (for increasing singularness of the integrand) increasingly
large local jumps --- that ruins the usefulness of $\hat{E}_4$, then $E_2$, and, for
non-integrable functions, finally even $E_1$.

From this excercise we conclude that it should always be a good idea, in any 
MC calculation, to monitor the behaviour of $E_2$ and $\hat{E}_4$ as $N$ increases;
and that this may tell us whether the second-order error, or indeed even the
first-order error itself, can be assigned any useful meaning. It should be pointed out
that, in our case study, the jumps in $\hat{E}_4$ are typically larger that those
in $E_2$ and that $\hat{E}_4$ is therefore a more sensitive probe of possible
convergence problems; and, independently of that, an estimate of how accurately the
integration error itself is estimated is in our opinion {\em always\/} adviseable.

\section*{Conclusions}
We have argued that the current practice of MC integration, resulting in a report
on the integral estimate and its error estimate, should always be accompanied by
a second-order error estimate, if only to validate the assignment of confidence
levels to the result (which can be, for instance, crucial in comparing the
results of different MC calculations, which is good and common practice). We have
presented the relevant estimators. A closer look at $E_4$ shows potential
positivitiy problems and we have emended this by defining an improved estimator
$\hat{E}_4$. We also point out that, on the one hand, the convergence of the
second-order error, ${\hat{E}_4}^{1/4}/{E_2}^{1/2} \sim N^{-1/4}$,
 rather than the `well-known' ${E_2}^{1/2}/E_1 \sim N^{-1/2}$
convergence of the error itself, and that on the other hand $E_2$ satisfies
its own version of the central-limit theorem. In addition, we have extended
the methods of the Chan-Golub-Leveque algorithm \cite{CG}
to allow for a numerically stable
computation of not only $E_2$ but $\hat{E}_4$ as well.

\section*{Appendix}
In this Appendix we will argue that the values of $\Sigma_2$ obey their own version
of  
CLT. This is not automatically obvious, since we can write
\bq
\Sigma_2 = N\suml_{j=1}^N\left(w_j-M(N)\right)^2
\eq
and therefore the summed quantities are not independent of one another.
Let us therefore consider a number of MC weights $w_j$, $j=1,2,\ldots,N$, that are 
identically distributed with probability density $P(w)$ but under the constraint
that 
\bq
\suml_{j=1}^Nw_j = 0\;\;.
\eq
We define 
\bq
X = {1\over N}\suml_{j=1}^N{w_j}^2\;\;,
\eq
and estimate the distribution of $X$ for large $N$ as follows. The moment-generating
function of $X$ reads
\bqa
\avg{e^{izX}} &\propto& \int du\;dw_1\cdots dw_N\;P(w_1)\cdots P(w_N)\;
\exp\left(iu\sum w_j + i{z\over N}\sum{w_j}^2\right)\nl
&=&\int du\left\{\int dw\;P(w)\;\exp\left(iuw + i{z\over N}w^2\right)\right\}^N\;\;,
\eqa
where the integrals run from $-\infty$ to $+\infty$. Introducing
\bq
\Phi_k(u) = \int dw\;P(w)\;e^{iuw}\;w^k
\eq
we can estimate
\bqa
\lefteqn{
\left\{\int dw\;P(w)\;\exp\left(iuw + i{z\over N}w^2\right)\right\}^N \;\;}\nl
&=& \exp\left(N\log\left(\Phi_0(u) + i{z\over N}\Phi_2(u)
-{z^2\over2N^2}\Phi_4(u) + {\cal O}\left({1\over N^3}\right)\right)\right)\nl
&\approx& \Phi_0(u)^N\;\exp\left(iz\lambda(u) - {z^2\over2N}\tau(u)\right)\;\;,\\
&&\lambda(u) = \Phi_2(u)/\Phi_0(u)\;\;\;,\;\;\;
\tau(u) = \Phi_4(u)/\Phi_0(u)\;\;.
\eqa
Now, since $\Phi_0(0)=1$ is the absolute maximum of $\Phi_0(u)$, and
\bq
\Phi_0(u) = 1 + iu\avg{w} - {u^2\over2}\avg{w^2} + {\cal O}(u^3)\;\;,
\eq
we can estimate
\bq
\left|\Phi_0(u)\right|^2 = 1 - u^2\sigma^2 + {\cal O}(u^4)\;\;\;,\;\;\;
\sigma^2 = \avg{w^2}-\avg{w}^2\;\;,
\eq
so that we may approximate
\bq
|\Phi_0(u)|^N \approx \exp\left(-{u^2N\over2}\sigma^2\right)
\eq
and the $u$ integral is dominated by the values of $u$ around zero; 
consequently,
\bqa
\avg{e^{izX}} \approx \exp\left(iz\lambda(0) - {z^2\over2N}\tau(0)\right)
\eqa
and the probability density for $X$ to take on the value $x$ is
\bqa
&&\mbox{Pr}(X=x) \;\;\propto\;\; \exp\left(-{N\over2\tau(0)}(x-\lambda(0))^2\right)\;\;,\nl
&&\lambda(0) = \avg{w^2}\;\;\;,\;\;\;\tau(0) = \avg{w^4}-\avg{w^2}^2\;\;.
\eqa
We see that in this sense a CLT holds for the distribution of $\sum{w_j}^2$.\\

As an illustration, we generate a large number ($10^6$) of samples of $N$ (pseudo-)random numbers uniformly in the interval $[0,1]$, and compute for these
$E_{1,2,4}$.  Below, we give the actual distribution of the $E_1$ values together with
the CLT Gaussian approximation with a width given by ${E_2}^{1/2}$. Similarly,
we also give the actual distribution of the $E_2$ values with their CLT Gaussian
approximation with width ${\hat{E}_4}^{1/2}$. We do this both for $N=10$ and for $N=1000$.

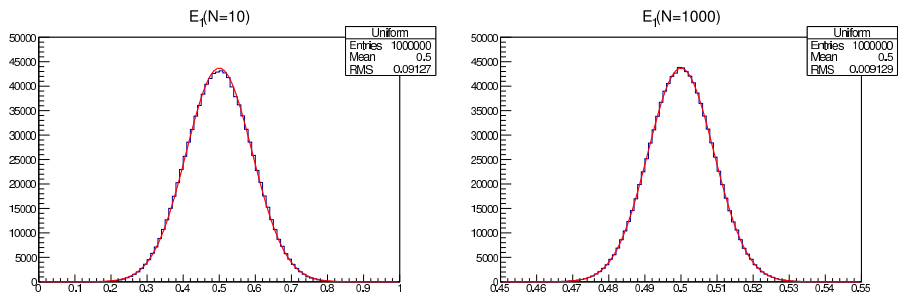
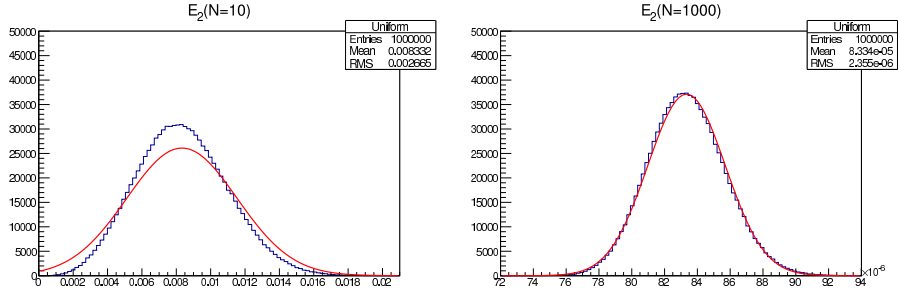
\begin{minipage}{6cm}{\[\diag{plotue110}{6}{0}\]}\end{minipage}
\begin{minipage}{6cm}{\[\diag{plotue11000}{6}{0}\]}\end{minipage}

\begin{minipage}{6cm}{\[\diag{plotue210}{6}{0}\]}\end{minipage}
\begin{minipage}{6cm}{\[\diag{plotue21000}{6}{0}\]}\end{minipage}

\noindent Unsurprisingly, for $N=1000$ the CLT approximation is excellent, but for
$N=10$ it is evdident that the approximation is much worse for $E_2$ than for 
$E_1$. 

We repeat the same excercise for numbers that are exponentially distributed,
that is, with probability density $P(x) = \exp(-x)$, $x\in[0,\infty)$. 

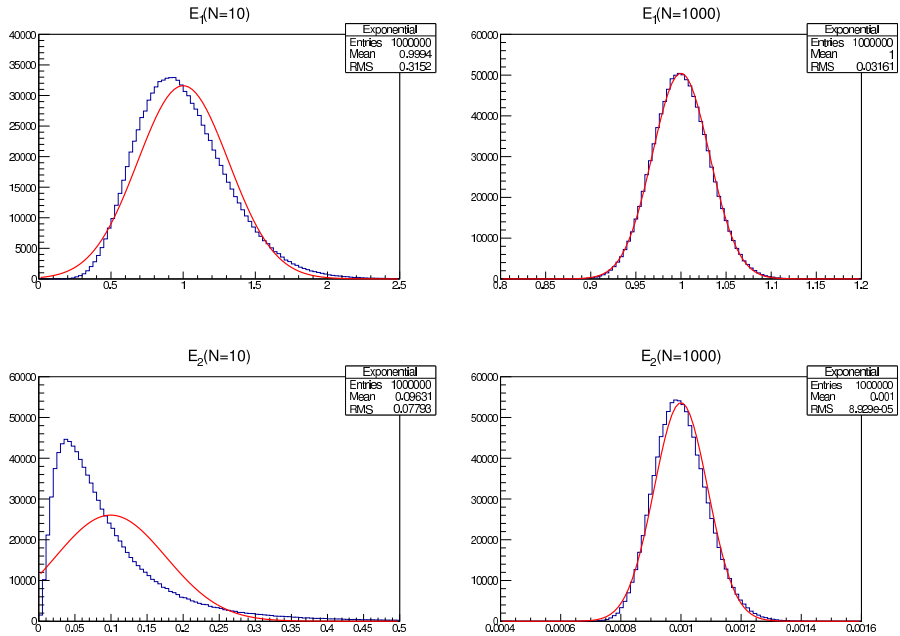
\begin{minipage}{6cm}{\[\diag{plotexe110}{6}{0}\]}\end{minipage}
\begin{minipage}{6cm}{\[\diag{plotexe11000}{6}{0}\]}\end{minipage}

\begin{minipage}{6cm}{\[\diag{plotexe210}{6}{0}\]}\end{minipage}
\begin{minipage}{6cm}{\[\diag{plotexe21000}{6}{0}\]}\end{minipage}

\noindent Because of the long high-$x$ tail of this $P(x)$, 
the CLT approximation is appreciably worse for $N=10$ although
still very good for $N=1000$.

Finally, we consider numbers distributed according to the exponential integral
\cite{abramowitz}:
\bq
P(x) = E_1(x) \equiv \intl_x^\infty dt\,{e^{-t}\over t}\;\;,\;\;x\in (0,\infty)\;\;,
\eq
which looks like $e^{-x}$ for large $x$, and like $-\log(x)$ for $x$ close to zero.
Such a distribution, with both many low-$x$ values and a high-$x$ tail, is
typical for how weights arising from MC event generators in particle physics
are distributed.

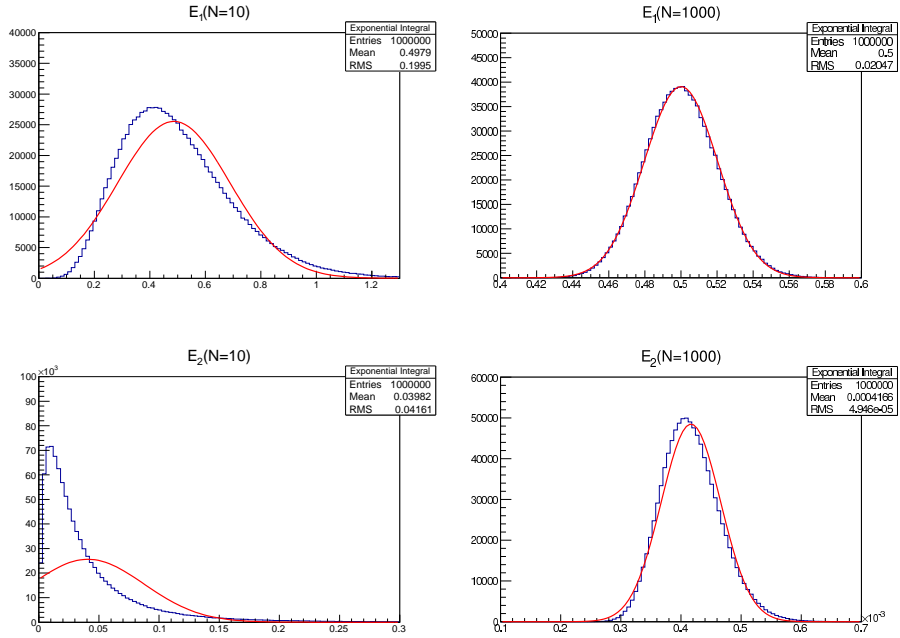
\begin{minipage}{6cm}{\[\diag{plote1e110}{6}{0}\]}\end{minipage}
\begin{minipage}{6cm}{\[\diag{plote1e11000}{6}{0}\]}\end{minipage}

\begin{minipage}{6cm}{\[\diag{plote1e210}{6}{0}\]}\end{minipage}
\begin{minipage}{6cm}{\[\diag{plote1e21000}{6}{0}\]}\end{minipage}

\noindent The CLT approximation is, unsurprisingly, very poor for  
$N=10$. However, for large $N$ values it is still seen to be quite good, where we 
must recall that $N=1000$ is actually quite a small number
for any serious calculation.


\begin{thebibliography}{9999}

\bibitem{MCoriginal} N. Metropolis, S. Ulam, {\it The Monte Carlo Method},
J. Am. Stat. Ass. 44 (1949) 335.

\bibitem{BKP}F.A.Berends, R.Pittau and R.Kleiss,
{\it Excalibur: A Monte Carlo program to evaluate all four fermion processes at LEP-200 and beyond},
    Comput.\ Phys.\ Commun.\  {\bf 85} (1995) 437.

\bibitem{James}  
F. James, {\it Monte Carlo Theory and Practice},
Rept.Prog.Phys. {\bf 43} (1980) 1145. 

\bibitem{CG}
%Tony F. Chan, Gene H. Golub, Randall J. Leveque,
T.F. Chan, G.H. Golub, R.J. Leveque,
{\it Algorithms for Computing the Sample Variance: Analysis and Recommendations,\/}
The American Statistician 37 (1983) 242.

\bibitem{abramowitz}
M. Abramowitz, and I.A. Stegun (eds),
{\it Handbook of Mathematical Functions\/}, Dover,1965.
%: with Formulas, Graphs, and Mathematical Tables


\end{thebibliography}
\end{document}